\documentclass[11pt,a4paper,reqno]{article} 
\usepackage{amsfonts,bm,pifont,alltt,amsthm,amscd,epsfig,amsmath,amssymb,enumerate,url,psfrag,txfonts,stmaryrd}

\newtheorem{theorem}{Theorem}
\newtheorem{definition}{Definition}
\newtheorem{proposition}[theorem]{Proposition}

\newtheorem{lemma}[theorem]{Lemma}

\def\card{\mathrm{Card}}

\def\eps{\varepsilon}
\def\II{\mbox{\rm  1\kern-0.20em I}}

\newcommand{\E}{\mathop{\hbox{\rm I\kern-0.17em E}}\nolimits}
\renewcommand{\P}{\mathop{\hbox{\rm I\kern-0.17em P}}\nolimits}

\newcommand{\Z}{\mathbb{Z}}
\newcommand{\Zb}{\overline{\mathbb{Z}}}
\newcommand{\N}{\mathbb{N}}
\newcommand{\R}{\mathbb{R}}

\newcommand{\Tc}{T_{\mathrm{col}}}

\begin{document} 

\title{Coexistence probability in the last passage percolation model is $6-8\log2$}

\author{\textsc{David Coupier}\footnote{\texttt{david.coupier@math.univ-lille1.fr}}\\
\textsc{Philippe Heinrich}\footnote{\texttt{philippe.heinrich@math.univ-lille1.fr}}}

\date{Laboratoire Paul Painlev\'e, UMR 8524\\
 UFR de Math\'ematiques, USTL, B\^at. M2\\
59655 Villeneuve d'Ascq Cedex France}

\maketitle

\noindent
{\bf Abstract:} A competition model on $\N^{2}$ between three clusters and governed by directed last passage percolation is considered. We prove that coexistence, i.e. the three clusters are simultaneously unbounded, occurs with probability $6-8\log2$. When this happens, we also prove that the central cluster almost surely has a positive density on $\N^{2}$. Our results rely on three couplings, allowing to link the competition interfaces (which represent the borderlines between the clusters) to some particles in the multi-TASEP, and on recent results about collision in the multi-TASEP.\\



\section{Introduction}

The directed \textit{last passage percolation} (LPP) model has been much
studied recently. In dimension $2$, it is closely related to some queueing
networks, to random matrix theory and to some combinatorial problems such as the longest increasing subsequence of a random permutation. See \textsc{Martin} \cite{Martin3} for a quite complete survey.

\begin{figure}[!ht]
\begin{center}
\begin{tabular}{cc}
\includegraphics[width=5.5cm,height=5.5cm]{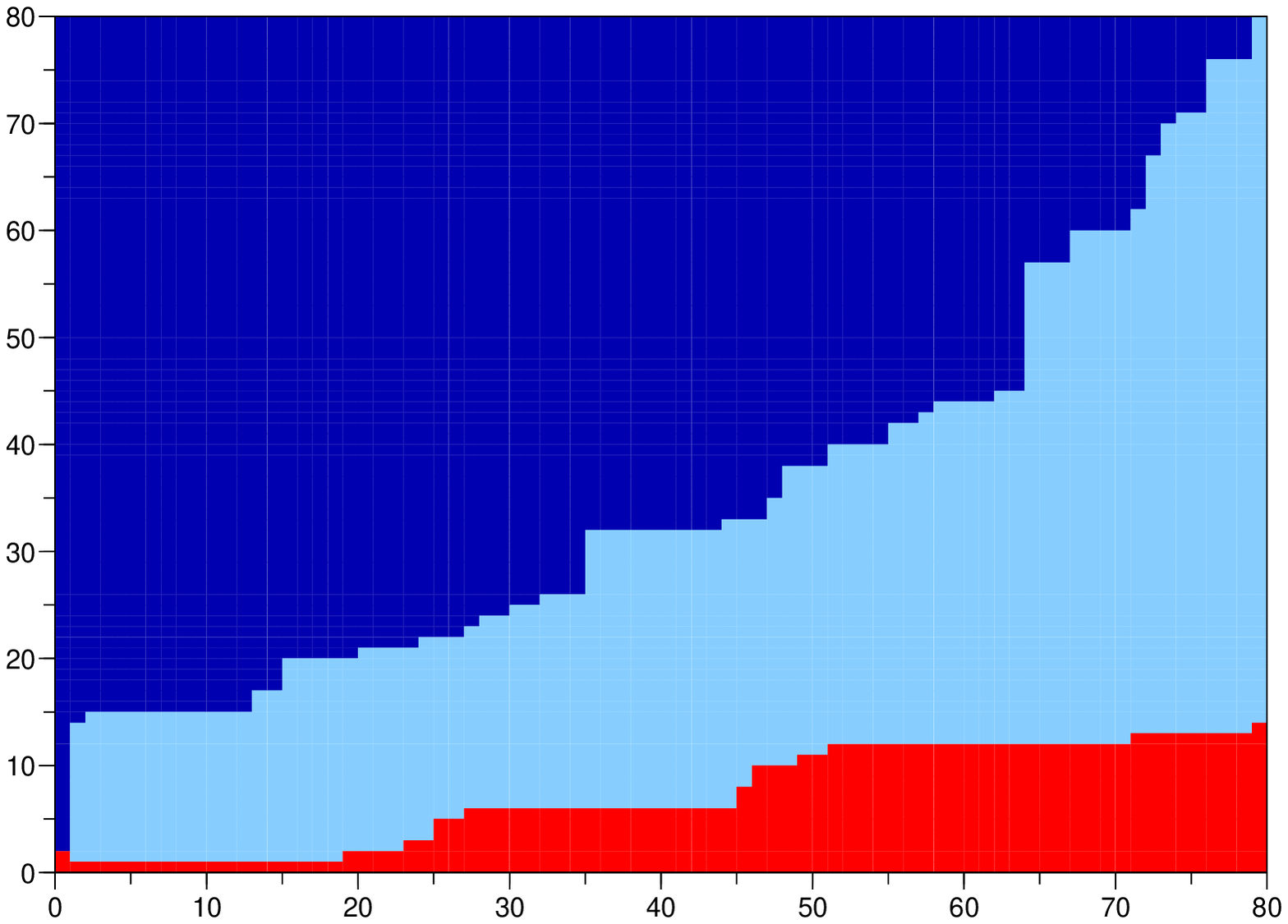} & \includegraphics[width=5.5cm,height=5.5cm]{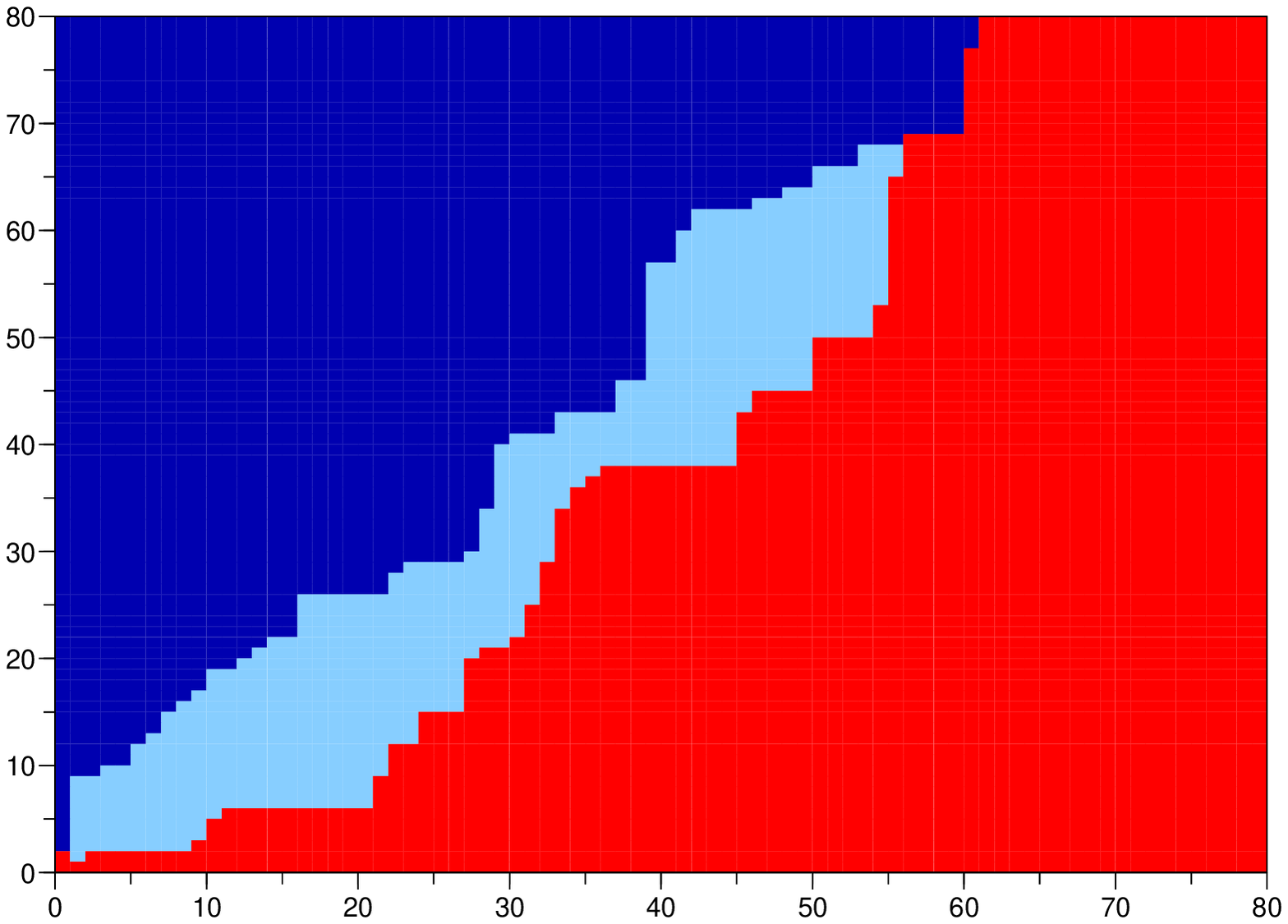}
\end{tabular}
\end{center}
\caption{\label{fig:simulations} Two simulations of the clusters $C(0,2)$,
  $C(1,1)$ and $C(2,0)$ which have been colored with respectively dark blue,
  light blue and red. To the left, $C(1,1)$ seems to be unbounded (there might
  be coexistence) whereas, to the right, it is bounded. Note that such a simulation of bounded but large cluster $C(1,1)$ is very rare.}
\end{figure}

Throughout this paper, $\N$ denotes the nonnegative integer set. We consider i.i.d. random variables $\omega(z)$, $z\in\N^{2}$, exponentially distributed with parameter $1$. Let $\P$ be the Borel probability measure induced by these variables on the product space $\Omega=[0,\infty)^{\N^{2}}$. The \textit{last passage time to $z$} is defined by
$$
G(z) = \max_{\gamma} \sum_{z'\in \gamma} \omega(z')
$$
where the above maximum is taken over all directed paths from the origin to
$z$ (see Section \ref{section:main} for precise definitions). The maximum
$G(z)$ is a.s. reached by only one path, called the \textit{geodesic} to
$z$. As a directed path, this geodesic goes through one and only one of the
three sites $(0,2)$, $(1,1)$ and $(2,0)$, called \textit{sources}. Let the
\emph{cluster} $C(s)$ be the set of sites $z\in\N^{2}$ whose geodesic goes by
the source $s$. Hence each configuration $\omega\in \Omega$ yields a random partition of $\{(x,y)\in\N^{2}:x+y\geq 2\}$, see Figure \ref{fig:simulations}.\\

\noindent
We focus on the competition (in space) between the three clusters $C(0,2)$, $C(1,1)$ and $C(2,0)$. The directed character of the model implies the first and the third ones are unbounded. But this is not necessary the case for the second one; we will talk about \textit{coexistence} when the cluster $C(1,1)$ is unbounded.\\
\indent
Our main result (Theorem \ref{theo:coex}) states that coexistence occurs with probability $6-8\log2$, which is close to $0.4548$. As far as we know, there is no other model where such a coexistence probability is exactly computed. For instance, in the (undirected) first passage percolation model, the competition between two clusters growing in the same space leads to two situations: either one cluster surrounds the other one, stops it and then infects all the other sites of $\mathbb{Z}^{2}$ or the two clusters grow mutually unboundedly, which is also called \textit{coexistence}. And in the
case of independent exponential weights, \textsc{H\"aggstr\"om} and \textsc{Pemantle} \cite{HP} have proved that coexistence occurs with positive probability. \textsc{Garet} and \textsc{Marchand} \cite{GM2} have since generalized this result to ergodic stationary passage times and to random environment.\\
Our second result (Theorem \ref{theo:strong}) completes the first one. When the cluster $C(1,1)$ is unbounded then it almost surely has a positive density in the following sense:
$$
\limsup_{n\to\infty} \frac{1}{n^{2}} \card \left( C(1,1) \cap [0,n]^{2} \right) > 0 ~,
$$
\indent
The proofs of Theorems \ref{theo:coex} and \ref{theo:strong} are mainly based
on three couplings: see \textsc{Thorisson} \cite{Th} for a complete reference
on couplings. The first one is due to \textsc{Rost}. In \cite{Rost}, he builds
a totally asymmetric simple exclusion process (TASEP) from the LPP model,
using the last passage times $G(z)$, $z\in\N^{2}$, as jump times. A background
on exlusion processes can be found in the book \cite{Liggett} (Part III) of \textsc{Liggett}. The borderlines between the clusters $C(1,1)$ and $C(0,2)$ and between $C(1,1)$ and $C(2,0)$ are modeled by two infinite directed paths, called the \textit{competition interfaces}. \textsc{Ferrari} and \textsc{Pimentel} \cite{FP}, thus \textsc{Ferrari}, \textsc{Martin} et al \cite{FMP} have studied their asymptotic behaviors. These competition interfaces play an important role here since the cluster $C(1,1)$ is bounded whenever they collide. The Rost's coupling allows to link these competition interfaces to two tagged pairs $[\infty \,1]$ in the TASEP, where labels $\infty$ and $1$ respectively represent holes and particles. In particular, the coexistence phenomenon is equivalent to the fact that these two tagged pairs never collide (Lemma \ref{lem:HIJ}).\\
The second coupling allows to turn the two tagged pairs into two \textit{second class particles} whose labels are denoted by $2$ and $3$ (Lemma \ref{lemma:TASEP/two-typeTASEP}). A second class particle is an extra particle which interacts with particles like a hole and interacts with holes like a particle. Its trajectory has been studied by \textsc{Mountford} and \textsc{Guiol} \cite{MountfordGuiol}. See also \textsc{Sepp{\"a}l{\"a}inen} \cite{Seppalainen2}. The idea to represent a second class particle as a hole-particle pair $[\infty\,1]$ is due to \textsc{Ferrari} and \textsc{Pimentel} \cite{FP}.\\
\textsc{Ferrari}, \textsc{Gon{\c{c}}alves} et al \cite{FGM} have studied the collision phenomenon of two second class particles. Thanks to the two previously announced couplings, they deduced (Theorem 4.1) that coexistence occurs in the LPP model with probability $1/3$. However, they assume for that some constraining initial conditions, namely
$\omega(0,0)=\omega(1,0)=\omega(0,1)=0$. We will explain why their coexistence result is a partial version of Theorem \ref{theo:coex}.\\
Finally, the third coupling, usually called basic coupling (\cite{Liggett},
p. 215), allows to consider the two second class particles
(i.e. the $2$ and $3$ particles) in a more general exclusion process, the
\textit{multi-TASEP}. Recently, \textsc{Amir} et al \cite{AAV} have proved
many results about this process. Some of them are expressed in terms of second
class particles (Proposition \ref{propo:collisiontim} and Lemma
\ref{lemm:Pm'}), thanks to that third coupling.\\
To sum up, these three couplings state a strong link between the multi-TASEP and the LPP model, leading to Theorems \ref{theo:coex} and \ref{theo:strong}.\\

The paper is organized as follows. Section \ref{section:main} contains the
definition of the LPP model and statements of  main results with some comments. Section \ref{Somedef} introduces the TASEP. The tagged pairs $[\infty \,1]$ are identified in Section \ref{taggedpairs}. Sections \ref{taggedthree-type} and \ref{multi-TASEP} are respectively devoted to the second and the third coupling on which proofs of Theorems \ref{theo:coex} and \ref{theo:strong} are based. The first coupling is described in Section \ref{RostCoupling}. Competition interfaces are defined in Section \ref{CompetInterf} and linked to tagged pairs in the TASEP in Section \ref{interftagged}. Finally, Theorems \ref{theo:coex} and \ref{theo:strong} are proved in Section \ref{proofs}.

\section{Coexistence results}
\label{section:main}

Recall that $\P$ denotes the law on $\Omega=[0,\infty)^{\N^{2}}$ of the family $\{\omega(z),z\in\N^{2}\}$ of i.i.d. random variables
exponentially distributed with parameter $1$.\\
A \textit{directed path} $\gamma$ from $(0,0)$ to $z$ is a finite sequence of
sites $(z_{0},z_{1},\ldots,z_{k})$ with $z_{0}=(0,0)$, $z_{k}=z$ and
$z_{i+1}-z_{i}=(1,0)\text{ or } (0,1)$, for $0\leq i\leq k-1$. The quantity
$\sum_{z'\in \gamma} \omega(z')$ represents the time to reach $z$ via $\gamma$. The set of all directed paths from $(0,0)$ to $z$ is denoted by $\Gamma(z)$. The \textit{last passage time to $z$} is defined by
$$
G(z) = \max_{\gamma\in\Gamma(z)} \sum_{z'\in \gamma} \omega(z') ~.
$$
Since each path of $\Gamma(z)$ goes by either $z-(1,0)$ or $z-(0,1)$, the function $G$ satisfies the recurrence relation
\begin{equation}
\label{recurrenceG}
G(z) = \omega(z) + \max \{ G(z-(1,0)) , G(z-(0,1)) \}
\end{equation}
(with boundary conditions $G(z)=0$ for $z=(x,-1)$ or $(-1,x)$ with
$x\in\N$). A site $z$ is said \textit{infected} at time $t$ if $G(z)\leq
t$. Relation (\ref{recurrenceG}) can be interpreted as follows: once both
sites $z-(1,0)$ and $z-(0,1)$ are infected, $z$ gets infected at rate $1$.

Recall that the cluster $C(s)$ is  the set of sites $z\in\N^{2}$ whose
geodesic goes by the source $s$. Let us point out the directed
character of the LPP model forces the clusters $C(2,0)$ and $C(0,2)$ to be
unbounded. Indeed, if the site $z=(x,y)$ belongs to $C(2,0)$,
so do all the sites on its right. Similarly, if the site $z=(x,y)$ belongs to
$C(0,2)$ so do all the sites above. Actually, only $C(1,1)$ can be bounded. Indeed, whenever
\begin{equation}
\label{coex<1}
\min \{ \omega(1,0)+\omega(2,0) , \omega(0,1)+\omega(0,2) \} > \omega(1,1) + \max \{ \omega(1,0) , \omega(0,1)\} ~,
\end{equation}
the last passage times $G(2,0)$ and $G(0,2)$ are both larger than $G(1,1)$. In this case, sites $(2,1)$ and $(1,2)$ respectively belong to $C(2,0)$ and $C(0,2)$, hence the cluster $C(1,1)$ is reduced to its source. See also the right hand side of Figure~\ref{fig:simulations} for the simulation of a larger (but bounded) cluster $C(1,1)$.\\
For any positive integer $n$, let
$$
\alpha(n) = \card \left( C(1,1) \cap \{(x,y) \in \N^{2}:\, x+y=n \} \right) ~.
$$
We will say there is \textit{coexistence} when the cluster $C(1,1)$ is
unbounded, i.e. $\alpha(n)>0$ for all $n\ge 2$. When this holds, each cluster $C(s)$ contains sites whose last passage time is as large as wanted; the three clusters $C(0,2)$, $C(1,1)$ and $C(2,0)$ coexist.

\begin{theorem}
\label{theo:coex}
Coexistence probability is $6-8\log2$:
$$
\P(\forall n\ge 2 , \; \alpha(n) > 0) = 6 - 8 \log 2 ~.
$$
\end{theorem}

It is already known that coexistence probability differs from $0$ and
$1$. Indeed, it is clear that coexistence cannot hold a.s. since the event
defined in (\ref{coex<1}) occurs with positive probability. Moreover, in a
previous work \cite{CH}, we have shown in particular that coexistence occurs
with positive probability if and only if there exists at least one infinite
geodesic (different from the horizontal and the vertical axes) with positive
probability; this last condition being proved in \cite{FP}, Proposition 7.

Let us compare our result to Theorem 4.1 of \cite{FGM}. In that paper, \textsc{Ferrari}, \textsc{Gon{\c{c}}alves} et al prove that coexistence occurs with probability $1/3$, but they consider the LPP model under the initial condition
\begin{equation}
\label{ci}
\omega(0,0) = \omega(1,0) = \omega(0,1) = 0 ~.
\end{equation}
Since the origin $(0,0)$ belongs to each geodesic, its weight does not affect
the coexistence probability. However, the cluster $C(1,1)$ benefits from
$\max\{\omega(1,0),\omega(0,1)\}$ whereas the clusters $C(2,0)$ and $C(0,2)$
only use $\omega(1,0)$ and $\omega(0,1)$ respectively. Assuming
$\omega(1,0)=\omega(0,1)=0$ amounts to remove this benefit. More precisely,
let $g:\Omega\to \Omega$ defined by $g(\omega)(0,0)=g(\omega)(1,0)=g(\omega)(0,1)=0$ and $g(\omega)(z)=\omega(z)$ otherwise. It then follows
$$
C(1,1)\left(g(\omega)\right) \subset C(1,1)(\omega) ~.
$$
Theorem 4.1 of \cite{FGM} says $C(1,1)\left(g(\omega)\right)$ is unbounded with probability $1/3$. This suggests that coexistence probability in the LPP model (without initial conditions) is greater than $1/3$. Actually, this remark has motivated the present work.\\

Our second result concerns the density of the cluster $C(1,1)$ in the quadrant
$\N^{2}$. Let us first remark that if the density of the cluster $C(1,1)$ is positive, i.e.
\begin{equation}
  \label{eq:positivedensity}
 \limsup_{n\to\infty} \frac{1}{n^{2}} \card \left( C(1,1) \cap [0,n]^{2} \right) > 0 ~, 
\end{equation}
then $C(1,1)$ is unbounded. Moreover, \eqref{eq:positivedensity} holds if and
only if 
$$\limsup_{n\to\infty} \frac{\alpha(n)}{n} > 0.$$
This stems from the fact that $\alpha(n+1)$ belongs to $\{\alpha(n)-1,\alpha(n),\alpha(n)+1\}$, for any $n$. Hence, the inequality $\alpha(n)>\delta n$ for some $\delta>0$ and integer $n$, implies that the cluster $C(1,1)\cap [0,n]^{2}$ contains a square with diagonal of length $\lfloor\delta n\rfloor$ (where $\lfloor x\rfloor$ denotes the integer part of $x$).
\begin{theorem}
\label{theo:strong}
Coexistence almost surely implies positive density for $C(1,1)$:
$$
\P \left.\left( \limsup_{n\to\infty} \frac{\alpha(n)}{n} > 0\; \right|\; \forall
  n\ge 2,\,\alpha(n) > 0 \right) = 1 ~.
$$
Moreover,
$$\P\left( \limsup_{n\to\infty} \frac{\alpha(n)}{n}<1\right) = 1 ~.$$
\end{theorem}

\section{TASEP and related processes}
\label{sectionTASEP}

\subsection{Some definitions}
\label{Somedef}

In the sequel, TASEP stands for totally asymmetric
simple exclusion process. It is a Markov process whose dynamics can be easily
described: at rate $1$ (i.e. after an exponential
time with parameter $1$), particles (integer or $\infty$) at sites $x$ and $x+1$ attempt to exchange their positions. The exchange occurs if the value at site $x$ is less than the value
at site $x+1$, otherwise nothing happens (total asymmetry property). There is at most one particle per site (exclusion condition). The $\infty$ particle has thus a role of hole. Here is a precise definition:

\begin{definition} Set $\Zb=\Z\cup\{\infty\}$. Let $S$ be a subset of $\Zb^\Z$. Consider the linear operator 
  $\mathcal{L}$ on cylinder functions $f$ on $S$ defined by
\begin{equation}
\label{generator}
\mathcal{L}f(\eta) = \sum_{x\in\Z} \II_{\{\eta_x<\eta_{x+1}\}} \left[ f\left(\eta^{x,x+1}\right) - f(\eta) \right] ~,
\end{equation}
where $\eta^{x,x+1}$ is obtained from $\eta=\{\eta_{y},y\in\Z\}$ by exchanging values at $x$ and $x+1$:
$$
\eta^{x,x+1}_{y} = \begin{cases}
\eta_{y} & \text{if $y\notin \{x,x+1\}$,}\\
\eta_{x+1} & \text{if $y=x$,}\\
\eta_{x} & \text{if $y=x+1$.}
\end{cases}
$$
A Markov process on $\R_+$ with configuration (or state) space $S$ and with generator
$\mathcal{L}$ is called 
\begin{enumerate}[(a)]
\item \emph{TASEP} if the configuration space  is  $S=\{1,\infty\}^\Z$,
\item \emph{$k$-type TASEP} if the configuration space is $S=\{1,2,\ldots,k,\infty\}^\Z$,
\item \emph{multi-TASEP} if the configuration space is $S=\Z^\Z$.
\end{enumerate}
\end{definition}

Let us add that the order relation $<$ on $\Z$ is extended to $\Zb$ as follows: $i<\infty$ if and only if $i$ belongs to $\Z$.\\
Besides, it will be convenient to locate some particles of interest in a configuration. Let $\eta$ be a configuration in $S\subset \Zb^\Z$ containing exactly one $k$ particle ($k\in\Zb$). The position of this $k$ particle in $\eta$ is denoted by
\begin{equation}
\label{keta}
k[\eta].
\end{equation} 
For a further use, it is convenient to introduce the following particular configurations, described in Figures \ref{fig:eta_m}, \ref{fig:eta_3m} and \ref{fig:basiccoupling}. For any integer $m$,
\begin{itemize}
\item Let $\eta^{m}\in \{1,\infty\}^\Z$ defined by
  \begin{equation}
    \label{eq:etam}
  \eta^{m}_{x} = \begin{cases}
1 & \text{ if $x\in \{\ldots,-3,-2\}\cup\{0\}\cup\{m+2\}$,}\\
\infty & \text{ otherwise.}
\end{cases}  
  \end{equation}
\item Let $\eta^{(3),m}\in \{1,2,3,\infty\}^\Z$ defined by
  \begin{equation}
\label{eq:eta3m}
  \eta^{(3),m}_{x} = \begin{cases}
1 & \text{ if $x\in \{\ldots,-3,-2,-1\}$},\\
2 & \text{ if $x=0$,}\\
3 & \text{ if $x=m+1$,}\\
\infty & \text{otherwise.}
\end{cases}  
  \end{equation}
\item Let $\eta^{(\infty)}\in \Z^\Z$ defined by
  \begin{equation}
    \label{eq:etainf}
  \eta^{(\infty)}_{x} = x \quad (x\in \Z) ~.
  \end{equation}
\end{itemize}

\subsection{Tagged pairs in the TASEP}
\label{taggedpairs}

\begin{figure}[!ht]
\begin{center}
\psfrag{m}{\small{$m$}}
\psfrag{-}{\small{$-$ pair}}
\psfrag{+}{\small{$+$ pair}}
\psfrag{1}{\small{$1$}}
\psfrag{i}{\small{$\infty$}}
\psfrag{Z}{\large{$\Z$}}
\includegraphics[width=9.2cm,height=1.7cm]{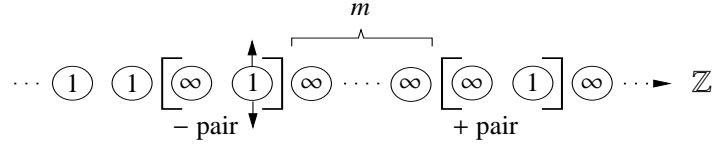}
\end{center}
\caption{\label{fig:eta_m} Configuration $\eta^{m}$ with the two tagged pairs $[\infty \,1]$. They are separated by $m$ ``holes'' $\infty$. On the axis $\Z$, the origin is marked with a vertical arrow.}
\end{figure}

We want to follow the evolution of two pairs of particles over time in the
TASEP with initial configuration $\eta^{m}$ defined in \eqref{eq:etam}. A pair
consists of a couple $(\infty,1)$ tagged with brackets. In the configuration
$\eta^{m}$, there are exactly two pairs $[\infty \,1]$, the left one is called
$-$ pair and the right one $+$ pair (see Figure \ref{fig:eta_m}).

Let us describe the evolution rule of the two pairs. Let $\eps\in\{-,+\}$. When a $1$ particle jumps (from the left and at rate $1$) over the hole $\infty$ of the $\eps$ pair, this one moves one unit to the left:
\begin{equation}
\label{evo1}
\overset{\curvearrowright}{1 [\infty} 1] \; \mbox{ becomes } \; [\infty\,1] 1~.
\end{equation}
When the $1$ particle of the $\eps$ pair jumps to the right (over a hole $\infty$ and at rate $1$) then the $\eps$ pair moves one unit to the right:
\begin{equation}
  \label{evo2}
[\infty \overset{\curvearrowright}{1]\infty} \; \mbox{ becomes } \; \infty [\infty\, 1] ~.
\end{equation}

\begin{definition}
\label{defi:Heps}
For $\eps\in\{-,+\}$, let us denote by $H^{\eps}(t)$ the hole's position of the $\eps$ pair, at time $t$, in
the TASEP with initial configuration $\eta^{m}$: 
$$
\cdots [\underset{H^-(t)}{\infty} 1] \cdots  [\underset{H^+(t)}{\infty} 1] \cdots
$$
The \emph{collision time} is defined as 
\begin{equation}
\label{Tc}
\Tc = \inf\{t\geq 0 : H^-(t) = H^+(t)\} ~,
\end{equation}
with the convention $\inf\emptyset=\infty$.
\end{definition}

The two tagged pairs merge together at the jump time $\Tc$ and there remains only one tagged pair so that $H^-(t)=H^+(t)=:H(t)$ for all $t\geq\Tc$:
$$
[\underset{H^-(t)}{\infty} 1\overset{\curvearrowright}{][}\underset{H^+(t)}{\infty}\,1] \; \mbox{ becomes } \;
\infty [\underset{H(t)}{\infty}\, 1] 1 ~.
$$
Finally, let us point out that the process $\{H^{\eps}(t),t\geq 0\}$, for $\eps\in\{-,+\}$, is not markovian but $\{(\xi(t),H^{-}(t),H^{+}(t)),t\geq 0\}$ is, where $\xi=\{\xi(t),t\geq 0\}$ is a TASEP with initial configuration
$\eta^{m}$. The reader can refer to \cite{Liggett} to get more details on tagged particle processes.

\subsection{From tagged pairs to three-type TASEP}
\label{taggedthree-type}

\begin{figure}[!ht]
\begin{center}
\psfrag{m}{\small{$m$}}
\psfrag{1}{\small{$1$}}
\psfrag{2}{\small{$2$}}
\psfrag{3}{\small{$3$}}
\psfrag{i}{\small{$\infty$}}
\psfrag{Z}{\large{$\Z$}}
\includegraphics[width=9.2cm,height=1.7cm]{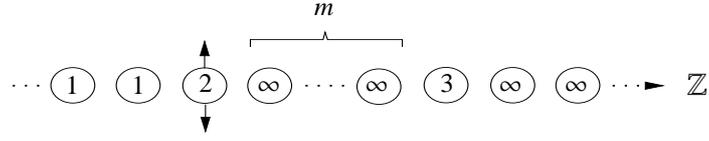}
\end{center}
\caption{\label{fig:eta_3m} Configuration $\eta^{(3),m}$. The two particles $2$ and $3$ are separated by $m$ holes $\infty$.}
\end{figure}

Recall that in a $k$-type TASEP (or in a multi-TASEP), a $i$ particle  can
pass a $j$ particle if and only if $i<j$. But the above evolution rule shows that
each tagged pair behaves like any single $i$ particle with respect to a $\infty$ particle --see \eqref{evo2}-- and also with respect to a  $1$ particle --see \eqref{evo1}-- provided $i$ is more than 1 and finite. If we turn the $-$ pair
into a $2$ particle and the $+$ pair
into a $3$ particle (for instance), we obtain a three-type TASEP. More
precisely, consider transformations
$$
\Psi^{x,y} = (\Psi^{x,y}_z)_{z\in\Z} : \{1,\infty\}^\Z \to \{1,2,3,\infty\}^\Z
$$
defined by 
\begin{enumerate}[(a)]
\item For $x+2\le y$, $$
\Psi^{x,y}_z(\eta) = \begin{cases}
\eta_{z-1} & \text{ if $z\leq x$,}\\
2 & \text{ if $z=x+1$,}\\
\eta_{z} & \text{ if $x+2\leq z\leq y-1$,}\\
3 & \text{ if $z=y$,}\\
\eta_{z+1} & \text{ if $z\geq y+1$.}
\end{cases}
$$
\item For $x=y$, $$
\Psi^{x,y}_z(\eta) = \begin{cases}
\eta_{z-1} & \text{ if $z\leq x-1$,}\\
3 & \text{ if $z=x$,}\\
2 & \text{ if $z=x+1$,}\\
\eta_{z+1} & \text{ if $z\geq x+2$.}
\end{cases}
$$
\end{enumerate}
For example,$\Psi^{-1,m+1}$ transforms $\eta^{m}$ (Figure \ref{fig:eta_m}) into $\eta^{(3),m}$ (Figure \ref{fig:eta_3m}). 

In what follows, we focus on the evolution of the two particles
$2$ and $3$  over time in the three-type TASEP until the collision time
$\Tc$. The applications $\Psi^{x,y}$ provide the following coupling:

\begin{lemma}
\label{lemma:TASEP/two-typeTASEP}
Let $\xi=\{\xi(t),t\geq 0\}$ be a TASEP with initial configuration
$\eta^{m}$ and collision time $\Tc$ as defined in \eqref{Tc}. Then, the process
$$
\xi':=\left\lbrace \Psi^{H^-(t),H^+(t)}(\xi(t)) , \; 0 \leq t \le \Tc\right\rbrace
$$
is a three-type TASEP on $[0,\Tc]$ with initial configuration
$\eta^{(3),m}$. In particular, with the notation \eqref{keta}, it follows
\begin{enumerate}[(i)]
\item For  $t<\Tc$,  $2[\xi'(t)]=H^{-}(t)+1$ and $3[\xi'(t)]=H^{+}(t)$,
\item For $t=\Tc$, $2[\xi'(t)]=H^{+}(t)+1=H^{-}(t)+1$ and  $3[\xi'(t)]=H^{+}(t)=H^{-}(t)$,
\end{enumerate}
\end{lemma}

It is crucial to remark this coupling holds until time $\Tc$ ($\Tc$ included thanks to the part (b) in the definition of $\Psi^{x,y}$).\\
Particles $2$ and $3$ in the three-type TASEP $\xi'$ can be seen as second class particles.

\subsection{From three-type TASEP to multi-TASEP}
\label{multi-TASEP}

\begin{figure}[!ht]
\begin{center}
\psfrag{Z}{\large{$\Z$}}
\psfrag{-2}{\small{$-2$}}
\psfrag{-1}{\small{$-1$}}
\psfrag{0}{\small{$0$}}
\psfrag{1}{\small{$1$}}
\psfrag{2}{\small{$2$}}
\psfrag{3}{\small{$3$}}
\psfrag{i}{\small{$\infty$}}
\psfrag{m}{\small{$m$}}
\psfrag{p}{\footnotesize{$m\negthickspace+\negthickspace2$}}
\psfrag{q}{\footnotesize{$m\negthickspace+\negthickspace1$}}
\includegraphics[width=10cm,height=3cm]{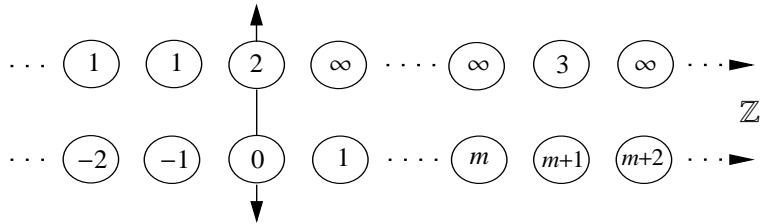}
\end{center}
\caption{\label{fig:basiccoupling} The configurations $\eta^{(3),m}$ and $\eta^{(\infty)}$ are the starting points of the three-type TASEP and the multi-TASEP under the basic coupling. They are aligned so that $2$ and $3$ particles in the three-type TASEP respectively correpond to $0$ and $m+1$ particles in the multi-TASEP (at time $t=0$).}
\end{figure}

The goal of this section is to couple a three-type TASEP with initial
configuration $\eta^{(3),m}$ and a multi-TASEP with initial
configuration $\eta^{(\infty)}$ using the \textit{basic coupling} (see \cite{Liggett}). To do so, let us consider a family $\{N_{x}(t), t\geq 0, x\in\Z\}$ of independent Poisson
processes with parameter $1$. At each event time $N_{x}(t)$ and for the two
processes, the particles located respectively at site $x$ and $x+1$ exchange their
positions if permitted by the order $<$, nothing changes
otherwise. Hence, the two processes evolve simultaneously on the
same probability space. See Figure \ref{fig:basiccoupling}.

First, let us remark some occurring jumps for the multi-TASEP, say between a $i$ particle and a $j$ particle (with $i<j$), are not authorized for the three-type TASEP. This happens when the corresponding particles in the three-type TASEP are the same or when $i\in\{1,\ldots,m\}$ and the corresponding particle to $j$ in the three-type TASEP is $3$. Then, we deduce that up to the time where the $2$ particle passes the $3$
one in the three-type TASEP,
\begin{itemize}
\item  the $2$ particle in the three-type TASEP
  corresponds to the $0$ particle  in the multi-TASEP;
\item  the $3$ particle in the three-type TASEP  corresponds to the further
  right particle among particles $1,\ldots,m+1$ in the multi-TASEP.
\end{itemize}
Hence, the  time where the $2$ particle and the $3$ particle exchange their
positions in the three-type TASEP is also the time where the $0$
particle has just overtaken all the particles $1,\ldots,m+1$ in the multi-TASEP. Theorem 7.1 of \textsc{Amir} et al. \cite{AAV} states this last event occurs with probability $2/(m+3)$. Now, the basic coupling allows to transfer this result to the three-type TASEP:

\begin{proposition}
\label{propo:collisiontim}
Let $\P_{m}'$ be the probability measure of a three-type TASEP $\xi'$ with initial configuration $\eta^{(3),m}$. With notation \eqref{keta}, it follows
$$
\P_{m}' \left(\exists t>0 , 2[\xi'(t)]>3[\xi'(t)]\right) = \frac{2}{m+3} ~.
$$
\end{proposition}

Note that, before results of \cite{AAV}, this result had been conjectured (and
proved in the case $m\in\{0,1\}$) by \textsc{Ferrari} et al in \cite{FGM}.\\

Let us respectively denote by $\xi'$ and $\xi''$ a three-type TASEP and a multi-TASEP with initial configurations $\eta^{(3),m}$ and $\eta^{(\infty)}$. Until the end of this section, we assume that $\forall t$, $2[\xi'(t)]<3[\xi'(t)]$. The basic coupling described above implies, at any time $t$, the $2$ particle in the three-type TASEP corresponds to the $0$ particle in the multi-TASEP, i.e.
$$
\forall t \geq 0 , \; 2[\xi'(t)] = 0[\xi''(t)] ~,
$$
and the $3$ particle in the three-type TASEP eventually corresponds to one of the particles $1,\ldots,m+1$ in the multi-TASEP, i.e.
$$
\exists k \in \{1,\ldots,m+1\} , \; \exists t_{k} , \; \forall t \geq t_{k} , \; 3[\xi'(t)] = k[\xi''(t)] ~.
$$
The fundamental result (Corollary 1.2) on which \cite{AAV} is based is that in the multi-TASEP with initial configuration $\eta^{(\infty)}$, each particle chooses a speed. Precisely, for every $k\in\Z$,
$$
\lim_{t\to\infty} \frac{k[\xi''(t)]}{t} = U_{k} \; \mbox{a.s.} ~,
$$
where $\{U_{k},k\in\Z\}$ is a family of random variables, each uniformly distributed on $[-1;1]$, called the \textit{TASEP speed process}. So, on the event $\{\forall t, 2[\xi'(t)]<3[\xi'(t)]\}$, the ratios $2[\xi'(t)]/t$ and $3[\xi'(t)]/t$ converge respectively to $U_{0}$ and $U_{k}$, for a given $k$. To sum up, the event
$$
\left\lbrace \lim_{t\to\infty} \frac{3[\xi'(t)]-2[\xi'(t)]}{t}=0 \; \mbox{ and } \; \forall t , 2[\xi'(t)]<3[\xi'(t)] \right\rbrace
$$
is a.s. included in
\begin{equation}
\label{equalspeed}
\bigcup_{k=1}^{m+1} \left\lbrace U_{0}=U_{k} \; \mbox{ and } \; \forall t , 0[\xi''(t)]<k[\xi''(t)] \right\rbrace ~.
\end{equation}
Finally, Lemma 9.9 of \cite{AAV} states, in the multi-TASEP with initial configuration $\eta^{(\infty)}$, every two particles with the same speed swap eventually. So the event (\ref{equalspeed}) has zero probability.

\begin{lemma}
\label{lemm:Pm'}
Let $\P_{m}'$ be the probability measure of a three-type TASEP $\xi'$ with initial configuration $\eta^{(3),m}$. Then,
$$
\P_{m}' \left( \lim_{t\to\infty} \frac{3[\xi'(t)]-2[\xi'(t)]}{t}=0 \; \mbox{ and } \; \forall t , 2[\xi'(t)]<3[\xi'(t)] \right) = 0 ~.
$$
\end{lemma}

\section{LPP model and tagged TASEP}
\label{sectionLPP}

The goal of this section is to state a coupling between the LPP model and the TASEP. This coupling allows to link the competition interfaces (defined in Section \ref{CompetInterf}) to some pairs of particles (identified in Section \ref{Initialconditions}).

\subsection{\textsc{Rost}'s coupling}
\label{RostCoupling}

\begin{figure}[!ht]
\begin{center}
\psfrag{1}{\small{$1$}}
\psfrag{i}{\small{$\infty$}}
\psfrag{c}{$P_{0}$}
\psfrag{b}{$P_{1}$}
\psfrag{a}{$P_{2}$}
\psfrag{d}{$H_{0}$}
\psfrag{e}{$H_{1}$}
\psfrag{f}{$H_{2}$}
\psfrag{Z}{\large{$\Z$}}
\includegraphics[width=8.5cm,height=2cm]{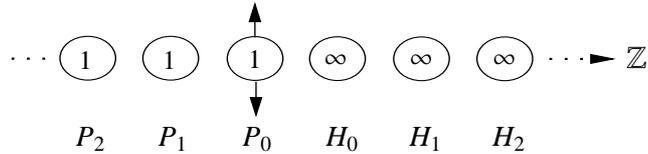}
\end{center}
\caption{\label{fig:rostinitial} Here is the configuration $\eta^{ext}$ with labelled particles. On the axis $\Z$, the origin is marked with the vertical arrow.}
\end{figure}

In \cite{Rost}, \textsc{Rost} gives an explicit construction of the TASEP from the LPP model, using the last passage times $G(z)$, $z\in\N^{2}$, as jump times. Let us describe this construction.\\
\indent
Let us start with the configuration $\eta^{ext}\in\{1,\infty\}^{\Z}$ which is made up of $1$ particles on nonpositive integers and $\infty$ particles on positive ones. The \textsc{Rost}'s idea consists in labelling $1$ particles from the right to the left by $P_{0}$, $P_{1}$, $P_{2}\ldots$ and $\infty$ particles from the left to the right by $H_{0}$, $H_{1}$, $H_{2}\ldots$ as in Figure \ref{fig:rostinitial} and in following them over time. Letters $P$ and $H$ refer to particle and hole.

\noindent
The evolution rule is
\begin{equation}
\label{Rostrule}
\mbox{$P_{j}$ and $H_{i}$ exchange their positions at time $G(i,j)$.}
\end{equation}
At time $G(0,0)=\omega(0,0)$ the first exchange takes place between $P_{0}$ and $H_{0}$. The second one will concern $P_{0}$ and $H_{1}$ if $\omega(1,0)<\omega(0,1)$ and $P_{1}$ and $H_{0}$ otherwise. More generally, at time $\max\{G(i-1,j),G(i,j-1)\}$, the exchanges between $P_{j}$ and $H_{i-1}$, and between $P_{j-1}$ and $H_{i}$ have already taken place. Labels of $1$ particles and those of $\infty$ particles remaining sorted over time then $P_{j}$ is then the left nearest neighbor of $H_{i}$. From that moment, they exchange their positions after the time $\omega(i,j)$ (i.e. at rate $1$) thanks to the recurrence relation (\ref{recurrenceG}):
$$
\omega(i,j) = G(i,j) - \max\{ G(i-1,j) , G(i,j-1) \} ~.
$$
It then suffices to disregard labels $P_{j}$ and $H_{i}$ to get back the TASEP. Precisely, let us denote by $P_{j}(t)$ and $H_{i}(t)$ the positions of particles $P_{j}$ and $H_{i}$ at time $t$. At the beginning, $P_{j}(0)=-j$ and $H_{i}(0)=i+1$. Now, set for $t\geq 0$ and $x\in\Z$
$$
\xi_x(t)=
\begin{cases}
1 & \mbox{if there exists $j$ such that $P_{j}(t)=x$,} \\
\infty & \mbox{otherwise,}
\end{cases}
$$
and let $\xi(t)$ be the configuration $(\xi_x(t))_{x\in\Z}$. Then:

\begin{lemma}
The process $\xi=\{\xi(t),t\geq0\}$ is the TASEP with initial configuration $\eta^{ext}$.
\end{lemma}

\begin{figure}[!ht]
\begin{center}
\psfrag{1}{\small{$1$}}
\psfrag{i}{\small{$\infty$}}
\psfrag{a}{$P_{0}$}
\psfrag{b}{$P_{1}$}
\psfrag{c}{$P_{2}$}
\psfrag{d}{$P_{3}$}
\psfrag{0}{$H_{0}$}
\psfrag{u}{$H_{1}$}
\psfrag{m}{$H_{m}$}
\psfrag{m1}{$H_{m\!+\!1}$}
\psfrag{m2}{$H_{m\!+\!2}$}
\psfrag{Z}{\large{$\Z$}}
\includegraphics[width=9.2cm,height=2cm]{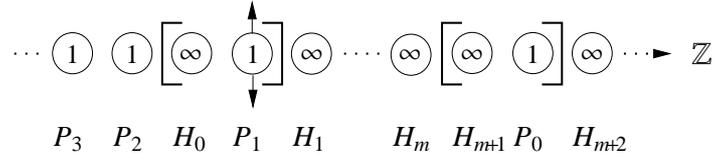}
\end{center}
\caption{\label{fig:rostm+1} On $\{N=m+1\}$ and at time $G(0,1)$, the TASEP obtained by the Rost's coupling is equal to the configuration $\eta^{m}$. The $-$ and $+$ pairs defined in Section \ref{taggedpairs} respectively consist of particles $H_{0}$ and $P_{1}$ and particles $H_{m+1}$ and $P_{0}$.}
\end{figure}

Let us end this section with describing an explicit way to obtain the configuration $\xi(t)$ from the infected region at time $t$, i.e. the set $\{z\in\N^{2}:G(z)\leq t\}$:
\begin{enumerate}
\item In the dual lattice $(-\frac{1}{2},-\frac{1}{2})+\N^{2}$, draw the border of the infected region at time $t$ and extend it on each side by two half-line, as in Figure \ref{fig:interfdebut}. The obtained broken line consists of  horizontal and vertical unit segments; it represents the axis $\mathbb{Z}$ on which $\xi(t)$ is defined.
\item Mark the last (from north to east) unit segment of the broken line before the diagonal $y=x$; it represents the origin of $\mathbb{Z}$.
\item Replace each vertical (resp. horizontal) unit segment of the broken line with a $1$ (resp. $\infty$) particle.
\end{enumerate}
For instance, the configuration of Figure \ref{fig:rostm+1} is obtained thanks to the previous algorithm from the infected region given by Figure \ref{fig:interfdebut}.

\subsection{Initial conditions in the LPP model}
\label{Initialconditions}

Consider the integer valued random variable $N$ defined by
$$
N = \left\{\begin{array}{ll}
\max \{ m \geq 1: \omega(1,0)+\ldots+\omega(m,0) < \omega(0,1) \} & \mbox{if exists,} \\
0 & \mbox{otherwise.}
\end{array} \right.
$$
We first remark that $\{N\ge 1\}=\{\omega(1,0)< \omega(0,1)\}$ occurs with probability  $1/2$, by symmetry.

\begin{lemma}
\label{lem:N}
Conditionally to $\{N\ge 1\}$, the random variable $N$ is distributed according to the geometric law with parameter $\frac{1}{2}$.
\end{lemma}

This result based on the memoryless property of the exponential distribution will be proved in Section \ref{prooflawN}.\\
Let $m$ be a nonnegative integer. The event $\{N=m+1\}$ implies that the first sites to be infected are in chronological order $(0,0),(1,0),\ldots,(m+1,0)$ and finally $(0,1)$; see Figure \ref{fig:interfdebut}. This provides the first
moves of particles in the TASEP $\xi$ obtained by the Rost's coupling. Precisely, $P_{0}$ overtakes $H_{0},\ldots,H_{m+1}$, thus at time $G(0,1)$ particle $P_{1}$ overtakes $H_{0}$. To sum up, on the event $\{N=m+1\}$, $\xi(G(0,1))$ is equal to the configuration $\eta^{m}$, introduced in (\ref{eq:etam}).

Since $G(0,1)$ is a stopping time, the strong Markov property implies

\begin{lemma}
\label{lem:PPm}
Conditionally to $\{N=m+1\}$, the shifted process $\xi(\cdot+G(0,1))$ is the TASEP with initial configuration $\eta^{m}$.
\end{lemma}

\subsection{Competition interfaces}
\label{CompetInterf}

Let us recall that $C(s)$ is the set of sites $z\in\N^{2}$ whose geodesic goes by the source $s$, for $s\in\{(0,2),(1,1),(2,0)\}$. The aim of this section is to define the borderlines between the clusters $C(2,0)$ and $C(1,1)$, and between $C(1,1)$ and $C(0,2)$.\\
\indent
The \textit{$-$ competition interface} is a sequence $(\varphi_{n}^{-})_{n\geq 0}$ defined inductively as follows: $\varphi_{0}^{-}=(0,0)$, $\varphi_{1}^{-}=(0,1)$ and for $n\geq 1$,
\begin{equation}
\label{varphi-}
\varphi_{n+1}^{-} = \left\{ \begin{array}{l}
\varphi_{n}^{-} + (1,0) \; \mbox{ if } \; \varphi_{n}^{-} + (1,1) \in C(0,2) ~,\\
\varphi_{n}^{-} + (0,1) \; \mbox{ if } \; \varphi_{n}^{-} + (1,1) \in C(1,1)\cup C(2,0) ~.
\end{array} \right.
\end{equation}
In an equivalent way, $\varphi_{n+1}^{-}$ chooses among the sites $\varphi_{n}^{-}+(1,0)$ and $\varphi_{n}^{-}+(0,1)$ the first to be infected. Moreover, it is easy to draw the competition interface $(\varphi_{n}^{-})_{n\geq 0}$ from a realization of clusters $C(2,0)$, $C(1,1)$ and $C(0,2)$. Indeed, $\varphi_{n}^{-}$ is the only site $(x,y)\in\N^{2}$ such that $x+y=n$, $(x+1,y)$ belongs to $C(1,1)\cup C(2,0)$ and $(x,y+1)$ to $C(0,2)$. So, the directed path $(\varphi_{n}^{-})_{n\geq 0}$ well describes the borderline between the clusters $C(1,1)$ and $C(0,2)$.\\
In the same spirit, the borderline between $C(2,0)$ and $C(1,1)$ is described by the \textit{$+$ competition interface}. This is a sequence $(\varphi_{n}^{+})_{n\geq 0}$ defined inductively by $\varphi_{0}^{+}=(0,0)$, $\varphi_{1}^{+}=(1,0)$ and for $n\geq 1$,
$$
\varphi_{n+1}^{+} = \left\{ \begin{array}{l}
\varphi_{n}^{+} + (1,0) \; \mbox{ if } \; \varphi_{n}^{+} + (1,1) \in C(1,1)\cup C(0,2) ~,\\
\varphi_{n}^{+} + (0,1) \; \mbox{ if } \; \varphi_{n}^{+} + (1,1) \in C(2,0) ~.
\end{array} \right.
$$
When the competition interfaces $(\varphi^{+}_{n})_{n\geq 0}$ and $(\varphi^{-}_{n})_{n\geq 0}$ meet on a given site $z_{0}$ (see the right hand side of Figure \ref{fig:simulations}) then they coincide beyond that site $z_{0}$ which is the larger (with respect to the $L^{1}$-norm) element of $C(1,1)$:
$$
\min \{ n \geq 1 , \; \varphi_{n}^{-} = \varphi_{n}^{+} \} = \max \{ x+y , \; (x,y) \in C(1,1) \} ~.
$$
In particular, there is coexistence if and only if the two competition interfaces never meet:
$$
\forall n \geq 2 , \; \varphi_{n}^{-} \not= \varphi_{n}^{+} ~.
$$

\subsection{From competition interfaces to tagged pairs}
\label{interftagged}

\begin{figure}[!ht]
\begin{center}
\psfrag{a}{\small{$0$}}
\psfrag{b}{\small{$1$}}
\psfrag{c}{\small{$2$}}
\psfrag{e}{\small{$m$}}
\psfrag{f}{\small{$m\!+\!1$}}
\includegraphics[width=7cm,height=3.5cm]{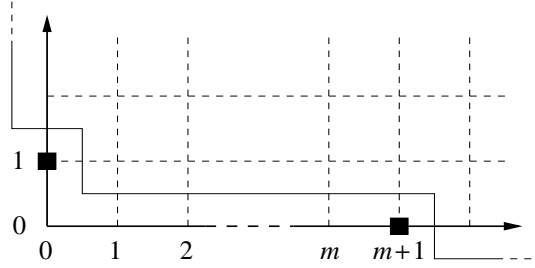}
\end{center}
\caption{\label{fig:interfdebut} The infected region at time $G(0,1)$ conditionally to $\{N=m+1\}$, deli\-mited by the black broken line. The two black squares represent $\phi^{-}(G(0,1))$ and $\phi^{+}(G(0,1))$. Since $\varphi_{n+1}^{\eps}$ chooses the earlier infected site among $\varphi_{n}^{\eps}+(1,0)$ and $\varphi_{n}^{\eps}+(0,1)$, the interface $\phi^{\eps}(t)$ is always in a corner formed by the black broken line. Combining with the algorithm given at the end of Section \ref{RostCoupling}, this justifies heuristically why $\phi^{\eps}(t)$ corresponds to a pair $[\infty\,1]$ in the TASEP $\xi$ obtained by the Rost's coupling.}
\end{figure}

Let $\eps\in \{+,-\}$. Consider the competition interface $(\varphi^{\eps}_{n})_{n\geq0}$ and its continuous-time counterpart, the interface process $\phi^{\eps}$ defined by
$$
\forall t \geq 0 , \; \phi^{\eps}(t) = \sum_{n\geq 0} \varphi^{\eps}_{n}
\mathbf{1}_{[G(\varphi^{\eps}_{n}),G(\varphi^{\eps}_{n+1}))}(t) ~.
$$
Set 
$$
\forall t \geq 0 , \; (I^{\eps}(t),J^{\eps}(t)):=\phi^{\eps}(t+G(0,1)) ~.
$$
By construction of $(\varphi^{-}_{n})_{n\geq0}$, $\phi^{-}(t)$ is $(0,0)$ until $G(0,1)$ and $\phi^{-}(G(0,1))$ is $(0,1)$. Besides, on the event $\{N=m+1\}$, the point $\phi^{+}(G(0,1))$ is known too. Assume this event satisfied. On the one hand, sites $(2,0),\ldots,(m+1,0)$ are infected before $(1,1),\ldots,(m,1)$ which yields $\varphi^{+}_{2}=(2,0),\ldots,\varphi^{+}_{m+1}=(m+1,0)$. On the other hand, at time $G(0,1)$, neither site $(m+2,0)$ nor site $(m+1,1)$ are still infected which means $\varphi^{+}_{m+2}$ is not yet determined. In conclusion, $\phi^{+}(G(0,1))$ is equal to $(m+1,0)$. To sum up, on the event $\{N=m+1\}$,
\begin{equation}
\label{I-t=0}
(I^{-}(0),J^{-}(0)) = (0,1) \; \mbox{ and } \; (I^{+}(0),J^{+}(0)) = (m+1,0) ~.
\end{equation}
See also Figure \ref{fig:interfdebut}.

Let $\xi$ be the TASEP obtained by the Rost's coupling and assume the event $\{N=m+1\}$ satisfied. Thanks to Lemma \ref{lem:PPm}, we know that the shifted process $\xi(\cdot+G(0,1))$ is the TASEP with initial configuration $\eta^{m}$. Recall that, in $\xi(t+G(0,1))$, the position of the $\infty$ particle of the $\eps$-pair is denoted by $H^{\eps}(t)$ (Definition \ref{defi:Heps}). Denote also by $P^{\eps}(t)$ the position of the $1$ particle of the $\eps$-pair. Of course, for any time $t$, $P^{\eps}(t)=H^{\eps}(t)+1$. Moreover, at time $t=0$ (and always on $\{N=m+1\}$),
\begin{equation}
\label{H-t=0}
(H^{-}(0),P^{-}(0))=(-1,0) \; \mbox{ and } \; (H^{+}(0),P^{+}(0))=(m+1,m+2) ~.
\end{equation}
The next result links competition interface $(\phi^{\eps}(t+G(0,1))_{t\geq 0}$ to the $\eps$ pair $[\infty\,1]$. Precisely, the coordinates $(I^{\eps}(t),J^{\eps}(t))$ are given by the labels of particles $\infty$ and $1$ constituting the $\eps$ pair at time $t$.

\begin{lemma}
\label{lem:HIJ}
The following identities hold on the event $\{N=m+1\}$. For any $t\geq0$ and $\eps\in\{+,-\}$,
\begin{equation}
\label{HIJ1}
\left( H^{\eps}(t) , P^{\eps}(t) \right) = \left( H_{I^{\eps}(t)}(t+G(0,1)) , P_{J^{\eps}(t)}(t+G(0,1)) \right) ~,
\end{equation}
and
\begin{equation}
\label{HIJ2}
H^{\eps}(t) = I^{\eps}(t)-J^{\eps}(t) ~.
\end{equation}
Moreover, for any $t\geq0$,
\begin{eqnarray}
\label{HIJ3}
H^{+}(t) = H^{-}(t) & \Leftrightarrow & \left( I^{+}(t) , J^{+}(t) \right) = \left( I^{-}(t) , J^{-}(t) \right) \\
& \Leftrightarrow & \phi^{+}(t+G(0,1)) = \phi^{-}(t+G(0,1)) \nonumber ~.
\end{eqnarray}
\end{lemma}

Recall that $\Tc$ is the time at which the tagged pairs collide (Definition \ref{defi:Heps}). Assume $N=m+1$ and $\Tc<\infty$. Then, just before time $\Tc$, the two tagged pairs in the TASEP $\xi(\cdot+G(0,1))$ are side by side and their labels satisfy $I^{-}(t)=I^{+}(t)-1$ and $J^{-}(t)=J^{+}(t)+1$ (thanks to (\ref{HIJ1})). Thus, at time $\Tc$,
the configuration $\cdots [\infty\,1][\infty\,1]\cdots$ becomes $\cdots\infty [\infty\,1] 1\cdots$ and thenceforward the two interfaces collide (thanks to (\ref{HIJ3})):
$$
\phi^{+}(\Tc+G(0,1)) = \phi^{-}(\Tc+G(0,1)) ~.
$$
Actually, $\Tc+G(0,1)$ is the time at which the last site of $C(1,1)$ is infected. Finally, remark the correspondence between competition interfaces and tagged pairs still holds after their collision.\\
Lemma \ref{lem:HIJ} will be proved in Section \ref{proofHIJ}.

\section{Proofs}
\label{proofs}

\subsection{Proof of Lemma \ref{lem:N}}
\label{prooflawN}

Let $m\geq 2$ be an integer. First,
$$
\P( N \geq m \;|\; N\ge 1 ) = \P( N \geq m \;|\; N \geq m-1 ) \times \P( N \geq m-1 \;|\; N\ge 1 ) ~.
$$
By the memoryless property of the exponential law,
\begin{eqnarray*}
\lefteqn{\P( N \geq m \;|\; N \geq m-1 )}\\
&\hspace{0.2cm} = & \P( \omega(1,0)+\ldots+\omega(m,0) < \omega(0,1) \;|\; \omega(1,0)+\ldots+\omega(m-1,0) < \omega(0,1) ) \\
&\hspace{0.2cm} = & \P( \omega(m,0) < \omega(0,1) ) \\
&\hspace{0.2cm} = & 1/2 ~.
\end{eqnarray*}
Hence, by induction we get $\P(N\geq m\;|\;N\ge 1)=2^{-m+1}$ which is also
true for $m=1$. This means that, conditionnally to  $\{N\ge 1\}$, $N$ is
geometrically distributed on $\{1,2,\ldots\}$ with parameter $1/2$. In other words,
$$
\P( N = m \;|\; N\ge 1 ) = 2^{-m} \quad (m\ge 1).
$$

\subsection{Proof of  Lemma \ref{lem:HIJ}}
\label{proofHIJ}

Throughout this proof, we assume $N=m+1$. Let us start with proving (\ref{HIJ1}) in the case $\eps=-$. In order to lighten formulas, let us denote by $\tau_{n}$ the time $G(\varphi^{-}_{n})-G(0,1)$. Since $\varphi^{-}_{1}=(0,1)$, $\tau_{1}$ is equal to $0$. At that time,
\begin{eqnarray*}
(H^{-}(0),P^{-}(0)) & = & (-1,0) \\
& = & (H_{0}(G(0,1)),P_{1}(G(0,1))) \\
& = & \left( H_{I^{-}(0)}(G(0,1)),P_{J^{-}(0)}(G(0,1)) \right) ~,
\end{eqnarray*}
thanks to relations (\ref{I-t=0}) and (\ref{H-t=0}). So, (\ref{HIJ1}) holds at time $\tau_{1}$ (and for $\eps=-$). Let us proceed by induction on times $(\tau_{n})_{n\geq1}$. Assume (\ref{HIJ1}) holds at time $\tau_{n}$ for a given $n\geq 1$, i.e. $I^{-}(\tau_{n})$ and $J^{-}(\tau_{n})$ are the labels of particles $\infty$ and $1$ of the $-$ pair at time $\tau_{n}$, and prove it still holds for any time $t\in[\tau_{n};\tau_{n+1}]$. By definition, $(I^{-}(\tau_{n}),J^{-}(\tau_{n}))$ are the coordinates of the competition interface $\phi^{-}(\tau_{n}+G(0,1))=\phi^{-}(G(\varphi^{-}_{n}))=\varphi^{-}_{n}$. At the next step, $\varphi^{-}_{n+1}$ chooses the earlier infected site among $(I^{-}(\tau_{n})+1,J^{-}(\tau_{n}))$ and $(I^{-}(\tau_{n}),J^{-}(\tau_{n})+1)$, say for example
$$
(I^{-}(\tau_{n+1}),J^{-}(\tau_{n+1})) = \phi^{-}(\tau_{n+1}+G(0,1)) = \varphi^{-}_{n+1} = (I^{-}(\tau_{n})+1,J^{-}(\tau_{n})) ~.
$$
Then, at time $\tau_{n+1}+G(0,1)=G(\varphi^{-}_{n+1})$, particles $H_{I^{-}(\tau_{n})+1}$ and $P_{J^{-}(\tau_{n})}$ exchange their positions while $H_{I^{-}(\tau_{n})}$ and $P_{J^{-}(\tau_{n})+1}$ have not yet done (see the Rost's rule (\ref{Rostrule})). This statement has two consequences. The first move of the $-$ pair after time $\tau_{n}+G(0,1)$ takes place at time $\tau_{n+1}+G(0,1)$: (\ref{HIJ1}) holds for any time $t\in[\tau_{n};\tau_{n+1})$. Thus, at time $\tau_{n+1}+G(0,1)$, the $-$ pair jumps one unit to the right and its particles $\infty$ and $1$ then become $H_{I^{-}(\tau_{n})+1}$ and $P_{J^{-}(\tau_{n})}$. So,
\begin{eqnarray*}
(H^{-}(\tau_{n+1}) , P^{-}(\tau_{n+1})) & = & (H_{I^{-}(\tau_{n})+1}(\tau_{n+1}+G(0,1)) , P_{J^{-}(\tau_{n})}(\tau_{n+1}+G(0,1))) \\
& = & (H_{I^{-}(\tau_{n+1})}(\tau_{n+1}+G(0,1)) , P_{J^{-}(\tau_{n+1})}(\tau_{n+1}+G(0,1))) ~,
\end{eqnarray*}
i.e. (\ref{HIJ1}) holds at time $\tau_{n+1}$. The case $\varphi^{-}_{n+1}=(I^{-}(\tau_{n}),J^{-}(\tau_{n})+1)$ leads to the same conclusion.\\
The case $\eps=+$ is very similar. This time, put $\tau_{n}=G(\varphi^{+}_{n})-G(0,1)$. We have already seen that at time $G(0,1)$ and on the event $\{N=m+1\}$, $\varphi^{+}_{m+1}=(m+1,0)$ and $\varphi^{+}_{m+2}$ is not yet determined. So $\tau_{m+1}<0$ and $\tau_{m+2}>0$. Relation (\ref{HIJ1}) holds at time $t=0$ thanks to (\ref{I-t=0}) and (\ref{H-t=0}):
\begin{eqnarray*}
(H^{+}(0),P^{+}(0)) & = & (m+1,m+2) \\
& = & (H_{m+1}(G(0,1)),P_{0}(G(0,1))) \\
& = & \left( H_{I^{+}(0)}(G(0,1)),P_{J^{+}(0)}(G(0,1)) \right) ~.
\end{eqnarray*}
Thus, the same induction as before, but on times $(0,\tau_{m+2},\tau_{m+3}\ldots)$, allows to conclude.\\
Let $\eps\in\{+,-\}$. It can be deduced from the previous remarks that when the $\eps$ pair jumps one unit to the right, i.e. $H^{\eps}$ increases by $1$, the label of its $\infty$ particle increases by $1$ whereas the one of its $1$ particle remains the same. Conversely, when the $\eps$ pair jumps one unit to the left, i.e. $H^{\eps}$ decreases by $1$, the label of its $1$ particle increases by $1$ whereas the one of its $\infty$ particle remains the same. To sum up, for any $t$,
$$
H^{\eps}(t)-H^{\eps}(0) = I^{\eps}(t)-J^{\eps}(t)-(I^{\eps}(0)-J^{\eps}(0)) ~.
$$
Combining with
$$
H^{-}(0) = -1 = I^{-}(0) - J^{-}(0) \; \mbox{ and } \; H^{+}(0) = m+1 = I^{+}(0) - J^{+}(0) ~,
$$
we get (\ref{HIJ2}).\\
It remains to prove (\ref{HIJ3}). Thanks to (\ref{HIJ2}), the equality $H^{+}(t)=H^{-}(t)$ is equivalent to
\begin{equation}
\label{I-I+}
I^{-}(t) - I^{+}(t) = J^{-}(t) - J^{+}(t) ~.
\end{equation}
Now, the directed character of the LPP model implies the differences $I^{-}(t)-I^{+}(t)$ and $J^{-}(t)-J^{+}(t)$ are respectively nonpositive and nonnegative. So, (\ref{I-I+}) forces $I^{-}(t)=I^{+}(t)$ and $J^{-}(t)=J^{+}(t)$.

\subsection{Proof of Theorem \ref{theo:coex}}
\label{section:proof1}

In Section \ref{CompetInterf}, the coexistence phenomenon has been described in terms of competition interfaces:
$$
\forall n \geq 2 , \; \alpha_{n} > 0 \; \Leftrightarrow \; \forall n \geq 2 , \; \varphi^{-}_{n} \not= \varphi^{+}_{n} ~.
$$
Let $m$ be a nonnegative integer. Relation (\ref{HIJ3}) of Lemma \ref{lem:HIJ} states that, on the event $\{N=m+1\}$, the two competition interfaces $(\varphi^{-}_{n})_{n\geq1}$ and $(\varphi^{+}_{n})_{n\geq1}$ never meet if and only if the collision time $\Tc$ of the tagged pairs in the shifted process $\xi(\cdot+G(0,1))$ obtained by the Rost's coupling, is infinite. Moreover, conditionally to $\{N=m+1\}$, $\xi(\cdot+G(0,1))$ is the TASEP with initial configuration $\eta^{m}$ (Lemma \ref{lem:PPm}). Let $\P_{m}$ be its probability measure. Then,
$$
\P \left( \forall n \geq 2 , \; \varphi^{-}_{n} \not= \varphi^{+}_{n} \;|\; N=m+1 \right) = \P_{m} (\Tc = \infty) ~.
$$
The coupling stated in Section \ref{taggedthree-type} between a TASEP with initial configuration $\eta^{m}$ and a three-type TASEP $\xi'$ with initial configuration $\eta^{(3),m}$ implies
$$
\P_{m} (\Tc = \infty) = \P_{m}' (\forall t , \; 2[\xi'(t)] < 3[\xi'(t)] ) ~,
$$
where $\P_{m}'$ denotes the probability measure of $\xi'$. Finally, the previous probability is equal to $1-2/(m+3)$ (Proposition \ref{propo:collisiontim}). Combining the previous identities, it follows:
\begin{equation}
\label{2/m+3}
\P \left( \forall n \geq 2 , \; \alpha_{n} > 0 \;|\; N=m+1 \right) = 1 - \frac{2}{m+3} ~.
\end{equation}
We conclude using symmetry of the LPP model, $\P(N\ge 1)=1/2$, Lemma \ref{lem:N} and (\ref{2/m+3}):
\begin{eqnarray*}
\P( \forall n \geq 2 , \; \alpha_{n} > 0 ) & = & 2 \P( \forall n \geq 2 , \; \alpha_{n} > 0 \;,\; N\ge 1 ) \\
& = & 2 \sum_{m=0}^{\infty} \P( \forall n \geq 2 , \; \alpha_{n} > 0 \;,\; N = m+1 ) \\
& = & 2 \sum_{m=0}^{\infty} \P( \forall n \geq 2 , \; \alpha_{n} > 0 \;|\; N = m+1 ) \\
& & \hspace*{3.5cm} \times \P( N = m+1 \;|\; N\ge 1 ) \P( N\ge 1 ) \\
& = & \sum_{m=0}^{\infty} \frac{1}{2^{m+1}} \left( 1 - \frac{2}{m+3} \right) \\
& = & 6 - 8 \log 2 ~.
\end{eqnarray*}
The last equality comes from the formula
$$
\log 2 = \sum_{m=1}^{\infty} \frac{1}{m 2^{m}} ~.
$$

Let us point out here that, thanks to the memoryless property of the exponential distribution, initial conditions $\omega(0,0)=\omega(1,0)=\omega(0,1)=0$ used in \cite{FGM} amounts to conditioning by the event $\{N=1\}$. So, their coexistence result (Theorem 4.1) corresponds to (\ref{2/m+3}) with $m=0$:
$$
\P( \forall n \geq 2 , \; \alpha_{n} > 0 \;|\; N = 1 ) = \frac{1}{3} ~.
$$

\subsection{Proof of Theorem \ref{theo:strong}}

Our goal is to prove that coexistence almost surely implies positive density:
\begin{equation}
\label{goalTh2}
\P \left( \forall n \geq 2 , \; \alpha(n) > 0 \; \mbox{ and } \; \lim_{n\to\infty} \frac{\alpha(n)}{n} = 0 \right) = 0 ~.
\end{equation}
For $\eps\in\{+,-\}$ and $n\geq 1$, let us denote by $\theta_{n}^{\eps}$ the angle formed by the half-line $[(0,0),\varphi_{n}^{\eps})$ with the axis $y=0$:
$$
\frac{\varphi_{n}^{\eps}}{|\varphi_{n}^{\eps}|} = e^{\mathbf{i} \theta_{n}^{\eps}} ~.
$$
Expressing the conditions $\alpha(n)>0$ and $\lim \alpha(n)/n=0$ in terms of angles $\theta^{-}_{n},\theta^{+}_{n}$ and using the symmetry of the LPP model with respect to the diagonal $x=y$, it is sufficient to prove
$$
\P \left.\left( \forall n \geq 2 , \; \theta^{-}_{n}>\theta^{+}_{n} \; \mbox{ and } \; \lim_{n\to\infty} \theta^{-}_{n}-\theta^{+}_{n} = 0 \;\right|\; N \geq 1 \right) = 0
$$
or, in an equivalent way, that the conditional probability
\begin{equation}
 \label{weakcoex}
\P \left.\left( \forall n \geq 2 , \; \theta^{-}_{n}>\theta^{+}_{n} \; \mbox{ and } \; \lim_{n\to\infty} \theta^{-}_{n}-\theta^{+}_{n} = 0 \;\right|\; N=m+1 \right)
\end{equation}
is null for any  $m\in\N$.\\
Let $m$ be a nonnegative integer. In \cite{FP}, \textsc{Ferrari} and \textsc{Pimentel} have studied the asymptotic behavior of the border between the two subsets $D(1,0)$ and $D(0,1)$ of $\N^{2}$ formed by sites whose geodesic respectively goes by $(1,0)$ and $(0,1)$. This border is described as a sequence $(\varphi_{n})_{n\geq0}$ --a competition interface-- defined by $\varphi_{0}=(0,0)$ and for $n\geq 0$,
$$
\varphi_{n+1} = \left\{ \begin{array}{l}
\varphi_{n} + (1,0) \; \mbox{ if } \; \varphi_{n} + (1,1) \in D(0,1) ~,\\
\varphi_{n} + (0,1) \; \mbox{ if } \; \varphi_{n} + (1,1) \in D(1,0) ~.
\end{array} \right.
$$
When $\omega(1,0)<\omega(0,1)$ the geodesic of $(1,1)$ goes by $(0,1)$ rather than $(1,0)$. In this case,
$$
D(1,0) = \{(1,0)\} \cup C(2,0) \; \mbox{ and } \; D(0,1) = \{(0,1)\} \cup C(0,2) \cup C(1,1) ~.
$$
So the sequences $(\varphi_{n})_{n\geq 0}$ and $(\varphi^{+}_{n})_{n\geq 0}$ coincide on the event $\{N=m+1\}$ which is included in $\{N\geq 1\}=\{\omega(1,0)<\omega(0,1)\}$. Now, by Proposition 4 of \cite{FP}, $(\theta^{+}_{n})_{n\geq 0}$ converges a.s. to a random angle $\theta$. Thus, by Proposition 5 of \cite{FP},
\begin{equation}
\label{I-J/t}
\lim_{t\to\infty} \frac{I^{+}(t)-J^{+}(t)}{t} = f(\theta) ~,
\end{equation}
where $f$ is a deterministic function (whose expression is without interest here). When the difference $\theta^{-}_{n}-\theta^{+}_{n}$ tends to $0$, results of \cite{FP} apply again and yield (\ref{I-J/t}) replacing $+$ with $-$. Therefore, (\ref{weakcoex}) is upperbounded by
$$
\P \left.\left(\begin{array}{c}
\lim_{t\to\infty} t^{-1} \left( I^{+}(t)-J^{+}(t) - (I^{-}(t)-J^{-}(t)) \right) = 0 \\
\mbox{ and } \; \forall t , \; I^{-}(t)-J^{-}(t) < I^{+}(t)-J^{+}(t)
\end{array} \;\right|\; N = m+1 \right) ~.
$$
Now, thanks to the Rost's coupling (Lemmas \ref{lem:PPm} and \ref{lem:HIJ}, relation (\ref{HIJ2})), the above conditional probability is equal to
\begin{equation}
\label{weakfinale}
\P_{m} \left( \lim_{t\to\infty} \frac{H^{+}(t)-H^{-}(t)}{t}=0 \; \mbox{ and } \; \forall t , \; H^{-}(t) < H^{+}(t) \right) ~,
\end{equation}
where $\P_{m}$ denotes the probability measure of the TASEP with initial configuration $\eta^{m}$. Finally, using Lemma \ref{lemma:TASEP/two-typeTASEP}, the quantity (\ref{weakfinale}) becomes
$$
\P_{m}' \left( \lim_{t\to\infty} \frac{3[\xi'(t)]-2[\xi'(t)]}{t}=0 \; \mbox{ and } \; \forall t , 2[\xi'(t)]<3[\xi'(t)] \right) ~,
$$
where $\xi'$ is a three-type TASEP with initial configuration $\eta^{(3),m}$ and $\P_{m}'$ its probability measure. Lemma \ref{lemm:Pm'} achieves the proof of (\ref{goalTh2}).\\
\indent
It remains to prove that a.s. the density of the cluster $C(1,1)$ cannot be equal to $1$. By symmetry with respect to the diagonal $x=y$, it suffices to show that
\begin{equation}
\label{density<1}
\P \left( \lim_{n\to\infty} \frac{\alpha(n)}{n} = 1 \; \mbox{ and } \; \omega(1,0) < \omega(0,1) \right) = 0 ~.
\end{equation}
When the density of the cluster $C(1,1)$ equals to $1$, that of cluster $C(2,0)$ is null. In this case, the $+$ competition interface $(\varphi^{+}_{n})_{n\geq0}$ is asymptotically horizontal and the sequence $(\theta^{+}_{n})_{n\geq0}$ converges to $0$. Furthermore, under the condition $\omega(1,0)<\omega(0,1)$, the competition interfaces $(\varphi^{+}_{n})_{n\geq0}$ and $(\varphi_{n})_{n\geq0}$ --previously introduced in this proof-- coincide. Proposition 4 of \cite{FP} then ensures the convergence almost sure of $(\theta^{+}_{n})_{n\geq0}$ to a random angle $\theta$. To sum up,
$$
\P \left( \lim_{n\to\infty} \frac{\alpha(n)}{n} = 1 \; \mbox{ and } \; \omega(1,0) < \omega(0,1) \right) \leq \P (\theta = 0) ~.
$$
Theorem 1 of \cite{FP} also says the distribution of $\theta$ has no atom. This proves (\ref{density<1}).\\

It derives from the above arguments that cluster $C(2,0)$ has a positive density on the event $\{\omega(1,0)<\omega(0,1)\}$, i.e. with probability one half. Actually, this holds with probability $1$ and the same is true for $C(0,2)$. To do so, let us remark that the cluster $C(2,0)$ grows when the weights $\omega(1,0)$ and $\omega(0,1)$ are exchanged, provided $\omega(1,0)$ is smaller than $\omega(0,1)$. It then can be proved that
$$
\P \left( \lim_{n\to\infty} \theta^{+}_{n} = 0 \; \mbox{ and } \; \omega(1,0) > \omega(0,1) \right) \leq \P \left( \lim_{n\to\infty} \theta^{+}_{n} = 0 \; \mbox{ and } \; \omega(1,0) < \omega(0,1) \right) ~.
$$
We have shown that the right hand side of the above inequality is null. Consequently, the probability of the event $\{\lim\theta^{+}_{n}=0\}$ is null which implies that the cluster $C(2,0)$ has a.s. a positive density.

\section*{Acknowledgments}

The authors thank P.~A. Ferrari and L.~P.~R. Pimentel for the communicating enthusiasm of their paper \cite{FP}. They also thank P.~A. Ferrari and J.~B. Martin for having focus our attention on some recent results in \cite{AAV}.

\bibliographystyle{plain}
\bibliography{LPP_bibli}

\end{document}